\documentclass[11pt]{article} 
 
\usepackage{latexsym} 
\usepackage{amsfonts} 
\usepackage{amsmath} 
\usepackage{amsthm} 
\usepackage{mathrsfs} 
\usepackage{dsfont} 
\usepackage{bbold} 
\usepackage[english]{babel} 
\usepackage{caption} 
\usepackage{epsfig} 
\usepackage{float} 
\usepackage{subfigure} 
\usepackage{psfrag} 
\usepackage{graphicx} 
\usepackage{epsfig} 
\usepackage{amsbsy} 
\usepackage{hyperref}
 
\DeclareMathOperator{\Val}{\matV}

\newtheorem{theorem}{Theorem} 
\newtheorem*{prop*}{Theorem} 
\newtheorem*{hyp*}{Hypothesis} 
\newtheorem{thm}{Theorem}

\newtheorem{lemma}[theorem]{Lemma}

\newtheorem{hyp}{Hypothesis}

\newcommand{\zerarcounters}{\setcounter{equation}{0}\setcounter{theorem}{0}}

\newcommand{\ZZZ}{\mathds{Z}} 
 
\newcommand{\NNN}{\mathds{N}} 
 
\newcommand{\RRR}{\mathds{R}} 
\newcommand{\TTT}{\mathds{T}}

\newcommand{\CCCC}{{\mathcal C}}

\newcommand{\calG}{{\mathcal G}}

\newcommand{\TT}{{\mathcal T}}

\newcommand{\gotn}{{\mathfrak n}}

\newcommand{\gotB}{{\mathfrak B}} 
\newcommand{\gotC}{{\mathfrak C}}

\newcommand{\matV}{{\mathscr V}}

\newcommand{\Fullbox}{{\rule{2.0mm}{2.0mm}}} 
 
\newcommand{\EP}{\hfill\Fullbox\vspace{0.2cm}} 
\newcommand{\prova}{\noindent{\it Proof. }} 
\newcommand{\io}{\infty} 
\newcommand{\e}{\varepsilon} 
\newcommand{\al}{\alpha} 
\newcommand{\de}{\delta} 
\newcommand{\be}{\beta}

\newcommand{\x}{\xi}

\newcommand{\oo}{\boldsymbol{\omega}}

\newcommand{\nn}{\boldsymbol{\nu}} 

\newcommand{\pps}{\boldsymbol{\psi}} 
\newcommand{\vzero}{\boldsymbol{0}}

\newcommand{\ii}{{\rm i}}

\def\ins#1#2#3{\vbox to0pt{\kern-#2 \hbox{\kern#1 #3}\vss}\nointerlineskip} 
 
\begin{document}
 
\title{\bf Response solutions for forced systems with large dissipation and arbitrary frequency vectors}
 
\author 
{\bf Guido Gentile and Faenia Vaia
\vspace{2mm} 
\\ \small
Dipartimento di Matematica, Universit\`a di Roma Tre, Roma, I-00146, Italy
\\ \small 
E-mail: gentile@mat.uniroma3.it
}

\date{} 
 
\maketitle 
 
\begin{abstract} 
We study the behaviour of one-dimensional strongly dissipative systems subject to
a quasi-periodic force. In particular we are interested in the existence of response solutions,
that is quasi-periodic solutions having the same frequency vector as the forcing term.
Earlier results available in the literature show that, when the dissipation is large enough
and a suitable function involving the forcing has a simple zero,
response solutions can be proved to exist and to be attractive provided
some Diophantine condition is assumed on the frequency vector.
In this paper we show that the results extend to the case of arbitrary frequency vectors.
 \end{abstract} 
  
\zerarcounters 
\section{Introduction} 
\label{sec:1} 

Consider the singular ordinary differential equation in $\RRR$
\begin{equation}\label{eq:1.1}
\e \ddot x + \dot x + \e \, g(x) = \e \, f(\oo t) ,
\end{equation}
where $\e\in\RRR$ is a small parameter and $\oo$ is a vector in $\RRR^{d}$, with $d\in\NNN$.
The functions $g:\RRR\to\RRR$ and $f:\TTT^{d}\to\RRR$ are assumed to be real analytic.
In particular, the function $t \mapsto f(\oo t)$ is quasi-periodic in $t$ and,
under the regularity assumptions we made, we can take its Fourier expansion
\begin{equation} \nonumber
f(\pps)=\sum_{\nn\in\ZZZ^{d}}{\rm e}^{\ii\nn\cdot\pps}f_{\nn} ,
\end{equation}
with the Fourier coefficients $f_{\nn}$ decaying exponentially in $\nn$.

For motivations and physical applications we refer to \cite{V,CLB,G1,CFG2} and references cited therein.
For a brief overview of related results on the existence of quasi-periodic solutions with frequency vectors
not satisfying any Diophantine condition we refer to Section \ref{sec:1.2} below.
Here we confine ourselves to recall that the equation \eqref{eq:1.1} with $\e>0$
describes a one-dimensional system in the presence of dissipation and
subject to an autonomous force $g$ and an additional quasi-periodic forcing term $f$.
The inverse of the \emph{perturbation parameter} $\e$ plays the role of the damping coefficient,
so that a small value for $\e$ corresponds to large dissipation -- on the contrary,
no smallness condition is assumed on the forces $f$ and $g$ acting on the system.
The vector $\oo$ is the \emph{frequency vector} of the forcing term. A \emph{response solution} to \eqref{eq:1.1}
is a quasi-periodic solution with the same frequency vector $\oo$ as the forcing.

\subsection{Main results} \label{sec:1.1}

If $\e=0$, for any constant $c\in \RRR$, $x=c$ is a solution to \eqref{eq:1.1}.
The problem we want to address is whether it is possible to choose the constant $c$ so that,
for $\e$ small enough, the equation \eqref{eq:1.1} admits a response solution
which tends to $c$ as $\e$ tends to zero.

Without any assumptions on the functions $f$ and $g$, the answer in general is negative.
More precisely, let us consider the following non-degeneracy condition.

\begin{hyp} \label{hyp:1}
Let $f$ and $g$ be the functions in \eqref{eq:1.1}.
The function $g(x)-f_{\vzero}$ has a simple zero $c_{0}$.
\end{hyp}

If the hypothesis is not satisfied one can provide counterexamples showing that response solutions
may fail to exist  \cite{G2}. On the contrary, if the hypothesis holds,
we shall prove that a response solution always exists, without requiring
any further condition on the frequency vector $\oo$.

This generalises previous results available in the literature and gives a positive answer
to a question raised in \cite{CFG2}. Indeed in \cite{GBD1,GBD2,G1,G2,CCD,CFG2}
some condition was assumed on $\oo$. We briefly recall the results.
If we set $\al_{n}(\oo) := \min \big\{ |\oo\cdot\nn| : 0<|\nn| \le 2^{n}  \big\}$ and define
\begin{equation} \nonumber
\e_{n}(\oo) := \frac{1}{2^{n}} \log \frac{1}{\al_{n}(\oo)} ,
\qquad \gotB(\oo) := \sum_{n=0}^{\io} \e_{n}(\oo) ,
\end{equation}
we say that $\oo\in\RRR^{d}$ satisfies
\begin{itemize}
\itemsep0em
\item the \emph{standard Diophantine condition} if
there exists two positive constants $\gamma_0$ and $\tau>0$ such that
$|\oo\cdot\nn| > \gamma_0 |\nn|^{-\tau}$ for all $\nn\in\ZZZ^{d}_{*} :=\ZZZ^{d}\setminus\{\vzero\}$,
\item the \emph{Bryuno condition} if the sequence $\e_{n}(\oo)$ is summable,
i.e. $\gotB(\oo)<\io$ \cite{B}.
\end{itemize} 
In \cite{GBD1} response solutions were proved to exist by assuming the vector $\oo$ to satisfy the
standard Diophantine condition and the result was then extended to vectors satisfying the weaker
Bryuno condition in \cite{GBD2} and to vectors $\oo$ such that $\e_{n}(\oo) \to 0$ as $n\to\io$ in \cite{CCD}.
The results was extended further to the case in which the zero $c_0$ of the function $g(x)-f_{\vzero}$
is of odd order $\gotn>1$, in \cite{G1,G2} for summable $\e_{n}(\oo)$ and in \cite{CFG2} for
$\e_{n}(\oo)$ converging to zero.

In this paper, still assuming Hypothesis \ref{hyp:1},
we remove the condition on the frequency vector and we only require that
the components of $\oo$ are rationally independent, that is
$\oo\cdot\nn\neq0$ $\forall\nn\in\ZZZ^{d}_{*}$.
This can be done without any loss of generality, since, if this is not the case,
$f$ can be expressed as a quasi-periodic function with frequency vector
$\oo'\in\RRR^{d'}$,  for some integer $d'<d$, with rationally independent components. 
Thus, we shall prove the following result.

\begin{thm} \label{thm:1}
Consider the ordinary differential equation \eqref{eq:1.1} and assume Hypothesis \ref{hyp:1}.
For any frequency vector $\oo\in\RRR^{d}$, there exists $\e_{0}>0$ such that for all $|\e|<\e_{0}$
there is at least one quasi-periodic solution $x_{0}(t)=c_{0}+X(\oo t,\e)$ to \eqref{eq:1.1},
such that $X(\pps,\e)$ is analytic in $\pps$ and goes to $0$ as $\e\to0$.
If $g'(c_0)>0$ the solution is locally attractive in the plane $(x,\dot x)$ and hence unique
in a neighbourhood of $(c_0,0)$.
\end{thm}

The equation \eqref{eq:1.1} is a special case of the more general
\begin{equation} \label{eq:1.2}
\e \ddot x + \dot x + \e \, h(x,\oo t) = 0 ,
\end{equation}
where $h:\RRR\times\TTT^{d}\to\RRR$ is a real analytic.
If we take the Fourier expansion of $\pps \mapsto h(x,\pps)$ by writing
\begin{equation} \nonumber
h(x,\pps)=\sum_{\nn\in\ZZZ^{d}}{\rm e}^{\ii\nn\cdot\pps} h_{\nn}(x) ,
\end{equation}
a natural counterpart of Hypothesis \ref{hyp:1} for \eqref{eq:1.2} is the following,

\begin{hyp} \label{hyp:2}
Let $h$ be the function in \eqref{eq:1.2}.
The function $h_{\vzero}(x)$ has a simple zero $c_{0}$.
\end{hyp}

Then the following result generalises Theorem \ref{thm:1}.

\begin{thm} \label{thm:2}
Consider the ordinary differential equation \eqref{eq:1.2} and assume Hypothesis \ref{hyp:2}.
For any frequency vector $\oo\in\RRR^{d}$, there exists $\e_{0}>0$ such that for all $|\e|<\e_{0}$
there is at least one quasi-periodic solution $x_{0}(t)=c_{0}+X(\oo t,\e)$ to \eqref{eq:1.2},
such that $X(\pps,\e)$ is analytic in $\pps$ and goes to $0$ as $\e\to0$.
\end{thm}

The two theorems provide information about the behaviour of a one-dimensional system
in the presence of large dissipation and subject to a quasi-periodic force,
as described by \eqref{eq:1.2} -- or by \eqref{eq:1.1} as a particular case.
The existence of a response solution depends on the zeroes of the function
$h_{\vzero}(x)$. If the function either has no zero or has a zero of even order,
the results in \cite{G2} show that in general no response solution exists.
On the contrary, if the function has a zero of order $\gotn=1$,
then the system always admits a response solution,
without any assumption on the frequency vector of the forcing term,
provided the dissipation is large enough. Such a result is obviously stronger
than the results of the papers quoted above, since no condition is assumed on $\oo$.
On the other hand, as it will emerge from the analysis of the next sections,
in general a smaller value is obtained for the estimate of $\e_{0}$:
this means that the closer $\oo$ to a resonance, the larger dissipation is needed
for the response solution to exist. Moreover, without any assumption on $\oo$,
less information is obtained about the regularity in $\e$ of the response solution. 
In fact, the stronger non-resonance condition on $\oo$, the more regularity is obtained on the dependence
of the solution on $\e$. For instance, if either $d=1$ (periodic case) or $d=2$ and $\oo$
satisfies a standard Diophantine condition with exponent $\tau=1$,
the solution turns out to be Borel-summable in $\e$ \cite{GBD2}.
Under the non-resonance condition considered in \cite{CFG2}, the response solution
is found to be $C^{\io}$ in $\e$ -- and analytic in a suitable domain with boundary tangent
to the origin \cite{CCD,CFG1}. On the contrary, if we do not assume any condition on $\oo$,
in general no more than a continuous dependence on $\e$ can be obtained.

We conclude with a few remarks.
\begin{enumerate}
\itemsep0em
\item
As mentioned above, if $\oo$ is such that the sequence $\e_{n}(\oo)$ converges to zero,
the existence of a response solution to \eqref{eq:1.1} can be proved under the weaker
assumption that order $\gotn$ of the zero of the function $g(x)-f_{\vzero}$ is odd \cite{CFG2}.
Unfortunately, the proof of theorem \ref{thm:1} given in Section \ref{sec:2} does not extend to the case $\gotn>1$.
\item \label{rmk2}
Another issue that deserves further investigation is the uniqueness and stability of the quasi-periodic solution.
In the case of Theorem \ref{thm:1}, under the assumption that $g'(c_0)>0$, asymptotic stability and hence
uniqueness follow from the same argument as given in \cite[Section 5]{BDG}.
\item
In this paper we have explicitly considered the one-dimensional case; however we expect the results to
carry over into the case in which $x\in\RRR^{n}$ and $f$ and $g$ are vector-valued real analytic function.
\end{enumerate}

The rest of the paper is organised as follows. We shall give the proof of Theorem \ref{thm:1} in Section \ref{sec:2}
and of Theorem \ref{thm:2} in Section \ref{sec:3}. Of course we could have confined ourselves to Theorem \ref{thm:2}, 
since Theorem \ref{thm:1} is included as a particular case. The reason why we have stated apart
Theorem \ref{thm:1} is that, as we shall see, the proofs in the forthcoming sections are in increasing order
of difficulty. Thus, it may be helpful to start considering first Theorem \ref{thm:1}, where the proof is easier,
before tackling the more general case where further technical intricacies arise. Moreover, the equation \eqref{eq:1.1}
studied in Theorem \ref{thm:1} is the one usually considered in the literature.

The proof of the theorems will be performed by introducing an auxiliary parameter $\mu$,
eventually to be put equal to 1, and looking for a formal power series expansion of the solution
in terms of $\mu$. Therefater, the series will be proved to be convergent by relying on a diagrammatic
representation of the coefficients and showing that the radius of convergence
is greater than 1 provide $\e$ is taken small enough: hence $\mu=1$ is allowed.
Note that the series is not a power series in $\e$. However, notwithstanding the solution,
as already observed above, is not even expected to be differentiable in $\e$,
we still are able to use convergent power series expansions.

\subsection{Related results on quasi-periodic solutions} \label{sec:1.2}

Periodic, quasi-periodic and almost periodic solutions in singularly perturbed systems have been
widely studied in the literature; we refer to \cite{BC2,Fi} for an introduction to almost periodic functions.
The ordinary differential equation \eqref{eq:1.1}, as well as \eqref{eq:1.2},
can be considered as a particular case of the system
\begin{equation} \label{eq:1.3}
\dot x = f(t,x,y,\e) , \qquad \e \, \dot y = g(t,x,y,\e) ,
\end{equation}
where, more generally, $(x,y) \in \RRR^{n}\times\RRR^{m}$, for some $n,m\in\NNN$,
and  the functions $f,g$ are almost periodic in time $t$. By assuming that an almost periodic
solution exists for the unperturbed system
\begin{equation} \label{eq:1.4}
\dot x = f(t,x,y,0) , \qquad 0 = g(t,x,y,0) ,
\end{equation}
then, under suitable non-degeneracy conditions on the vector field,  an almost periodic is shown to exist nearby;
see for instance \cite{FL,A,W,V,F} in the periodic case and \cite{HS,C,S} in the almost periodic case.
Note that the non-degeneracy conditions are not satisfied by our equation.
Indeed, while the quoted results can be seen as results on the persistence of central manifolds under
suitable hypothesis of stability or hyperbolicity or exponential dichotomy
(for instance attractive limit cycles in van de Pol-like systems),
theorems \ref{thm:1} and \ref{thm:2} are about the bifurcation of quasi-periodic solutions from a fixed point.
 
Almost periodic solutions were also studied in non-singularly perturbed linear and nonlinear systems, once more
assuming exponential dichotomy on the unperturbed linearised system; see for instance  \cite{H,C,F,He,XLC}.

In the different context of conservative systems, nonlinear Duffing equations
with almost periodic forcing were studied in \cite{BC1,BC2,BZ} with variational techniques.
In particular, almost periodic solutions with the same frequency vector as the forcing were proved to exist
without assuming any smallness condition on the forcing and any Diophantine condition on its frequency vector.
The latter issue marks a remarkable difference with respect KAM-like results, where restrictive
conditions are imposed on the frequency vector. Note, however, then,
even in the perturbative regime (small forcing), the solutions do not describe KAM invariant curves,
as they bifurcate from a stable fixed point which they reduce to in the absence of the forcing --
a situation which has some analogies with the kind of solutions we discuss in this paper.

In all the aforementioned papers, uniqueness of the the solution is proved as well.
However, the hypotheses on the equations assumed in the papers ensure that the solutions
are bounded for all times -- a property that in the case of equation \eqref{eq:1.1}
holds if $g'(c_0)>0$. Under such further assumption, as noted in the second remark at the end
of Section \ref{sec:1.1}, also the quasi-periodic solution to \eqref{eq:1.1}
turns out to be unique (in fact locally attractive).

\zerarcounters 
\section{Proof of Theorem \ref{thm:1}}
\label{sec:2} 

Henceforth, we assume $\oo$ to be non-resonant, that is $\oo\cdot\nn \neq 0$
for all $\nn\in\ZZZ^{d}_{*}=\ZZZ^{d}\setminus\{\vzero\}$.
As noted in Section \ref{sec:1} this is not restrictive.
We shall prove that there exists a response solution $x_0(t)$ to \eqref{eq:1.1}
such that $(x_0(t),\dot x(t))$ describes a curve in a neighbourhood of $(c_0,0)$.
As already observed in the second remark at the end of Section \ref{sec:1.1},
stability and uniqueness can be proved exactly as in \cite{BDG}.

Let us denote by $\Sigma_{\xi}$ the strip of $\TTT^{d}$ of width $\xi$ and
by $\Delta(c_{0},\rho)$ the disk of center $c_{0}$ and radius $\rho$ in the complex plane. 
By the assumptions on $f$ and $g$, for any $c_{0}\in\RRR$ there exist
$\xi_{0}>0$ and $\rho_{0}>0$ such that $\pps \mapsto f(\pps)$ is analytic
in $\Sigma_{\xi_{0}}$ and $x \mapsto g(x)$ is analytic in $\Delta(c_{0},\rho_{0})$.
Then for all $\xi<\xi_{0}$ and all $\rho<\rho_{0}$ one has
\begin{subequations} \label{eq:2.1}
\begin{align}
f(\pps) & = \sum_{\nn\in\ZZZ^{d}} {\rm e}^{\ii\nn\cdot\pps} f_{\nn} , \qquad 
|f_{\nn}| \le \Phi \, {\rm e}^{-\x |\nn|} ,
\label{eq:2.1a} \\
g(x) & = g(c_0) + \sum_{p=\gotn}^{\io} a_{p} (x-c_{0})^{p} , \qquad 
a_{p} := \frac{1}{p!} \frac{{\rm d}^{p}g}{{\rm d}x^{p}}(c_{0}) , \qquad
|a_{p}| \le \Gamma \, \rho^{-p} ,
\label{eq:2.1b}
\end{align}
\end{subequations}
where $\Phi$ is the maximum of $f(\pps)$ for $\pps\in\Sigma_{\xi}$
and $\Gamma$ is the maximum of $g(x)$ for $x\in\Delta(c_{0},\rho)$.

Let us rewrite (\ref{eq:1.1}) as
\begin{equation} \label{eq:2.2}
\e \ddot x + \dot x + \e \, g(c_0) + \e \, a \, (x-c_{0}) + \e \, G(x) = \e \, f(\oo t) ,
\end{equation}
where $a:=a_{1} \neq0$ by Hypothesis \ref{hyp:1}, and
\begin{equation} \nonumber 
G(x) := g(x)-g(c_{0})-a\left( x - c_0 \right) = \sum_{p=2}^{\io} a_{p} (x-c_0)^{p} , \\
\end{equation}

We look for a quasi-periodic solution to (\ref{eq:2.2}), that is a solution of the form
\begin{equation} \label{eq:2.3}
x(t,\e) = c_{0} + \zeta + u(\oo t, \e,\zeta) , \qquad u(\pps,\e,\zeta) = 
\sum_{\nn\in\ZZZ^{d}_{*}} {\rm e}^{\ii\nn\cdot\pps} u_{\nn} ,
\end{equation}
where $\zeta$ is a parameter that has to be fixed eventually and $\pps \mapsto u(\pps,\e,\zeta)$
is a zero-average quasi-periodic function, with Fourier coefficients depending on both $\e$ and $\zeta$.
Thus, we can write \eqref{eq:2.2} in Fourier space. If we set
\begin{equation} \label{eq:2.4}
D(\e,s) := - \e s^{2} + \ii s + \e \, a ,
\end{equation}
we obtain the following equations:
\begin{subequations} \label{eq:2.5}
\begin{align}
D(\e,\oo\cdot\nn) \, u_{\nn} & = - \e \left[ G(c_0+\zeta+u) \right]_{\nn} + \e \, f_{\nn} , \qquad \nn \neq \vzero ,
\label{eq:2.5a} \\
\e \, a \, \zeta  & = - \e \left[ G(c_0+\zeta+u) \right]_{\vzero} ,
\label{eq:2.5b}
\end{align}
\end{subequations}
where we have used that $g(c_0)=f_{\vzero}$ by Hypothesis \ref{hyp:1}.
By the notation $[G(c_0+\zeta+u)]_{\nn}$ we mean that
we first write $u$ according to \eqref{eq:2.3}, then expand $G(c_0+\zeta+u)$ in Fourier series
in $\pps$ and finally keep the Fourier coefficient with index $\nn$.
The splitting into two sets of equations is typical of the Lyapunov-Schmidt reduction \cite{CH,GS,SLSF};
we call \eqref{eq:2.5a} the \emph{range equation} and \eqref{eq:2.5b} the \emph{bifurcation equation}.

In order to solve \eqref{eq:2.5}, we proceed as follows. We first ignore \eqref{eq:2.5b} and look for a solution
to \eqref{eq:2.5a}, depending on the parameter $\zeta$. If we are able to do this, then
we pass to \eqref{eq:2.5b} and try to fix $\zeta$ in such a way to make such an equation to be satisfied.

\subsection{The range equation}

We start by studying \eqref{eq:2.5a} only and considering $\zeta$ as a free parameter,
close enough to $0$. We introduce the auxiliary parameter $\mu$ by modifying (\ref{eq:2.5a}) into
\begin{equation} \label{eq:2.6}
D(\e,\oo\cdot\nn) \, u_{\nn} = - \mu \left[ \e G(c_0+\mu\zeta + u) \right]_{\nn} + 
\mu \e \, f_{\nn} , \qquad \nn \neq \vzero ,
\end{equation}
and look for a quasi-periodic solution to \eqref{eq:2.6} in the form of a power series in $\mu$,
\begin{equation} \label{eq:2.7}
u(\oo t, \e,\zeta,\mu) = \sum_{k=1}^{\io} \sum_{\nn\in\ZZZ^{d}_{*}} 
\mu^{k} {\rm e}^{\ii\nn\cdot\pps} u^{(k)}_{\nn} .
\end{equation}
We shall show that there exists $\mu_{0}>0$ such that, for all $\zeta$ small enough,
there exists a solution of the form \eqref{eq:2.7}, analytic in $\mu$ for $|\mu|<\mu_{0}$.
The original equation \eqref{eq:2.5a} is recovered when $\mu=1$, so we need $\mu_{0}>1$. 
The argument below is a variant of that given in \cite{CFG2}, to which we refer for more details
about the construction described hereafter.

By inserting \eqref{eq:2.7} into \eqref{eq:2.6} we obtain a recursive definition
for the coefficients $u^{(k)}_{\nn}$. To simplify the notations, we set $u^{(1)}_{\vzero}=\zeta$
and $u^{(k)}_{\vzero}=0$ $\forall k\ge 2$. Then, one has, formally,
\begin{subequations} \label{eq:2.8}
\begin{align}
D(\e, \oo\cdot\nn) \, u^{(1)}_{\nn} & = \e \, f_{\nn} 
\label{eq:2.8a} \\
D(\e, \oo\cdot\nn) \, u^{(k)}_{\nn} & = - \e \, \sum_{p=2}^{\io} a_{p} \!\!\!\!\!\!
\sum_{\substack{ k_{1},\ldots,k_{p} \ge 1 \\ k_{1}+\ldots+k_{p}=k-1}}
\sum_{\substack{\nn_{1},\ldots,\nn_{p} \in\ZZZ^{d} \\ \nn_{1}+\ldots+\nn_{p}=\nn}} \!\!\!\!\!\!
u^{(k_{1})}_{\nn_{1}} \ldots u^{(k_{p})}_{\nn_{p}} , \qquad k \ge 2 ,
\label{eq:2.8b}
\end{align}
\end{subequations}
where we recall that $\nn\neq \vzero$. Here and henceforth the sums over the empty set are meant as zero.
In particular (see Remark 2.1 in \cite{CFG2})  for $k=2$ one has $u^{(2)}_{\nn}=0$ $\forall\nn\in\ZZZ^{d}_{*}$.

By iterating \eqref{eq:2.8b} one obtains a diagrammatic representation of the coefficients
$u^{(k)}_{\nn}$ in terms of trees. The construction is very similar to that in \cite{CFG2};
see also \cite{G3} for a review on the tree formalism.

A \textit{rooted tree} $\theta$ is a graph with no cycle,
such that all the lines are oriented toward a unique
point (\textit{root}) which has only one incident line (root line).
All the points in $\theta$ except the root are called \textit{nodes}.
The orientation of the lines in $\theta$ induces a partial ordering 
relation ($\preceq$) between the nodes. Given two nodes $v$ and $w$,
we shall write $w \prec v$ every time $v$ is along the path
(of lines) which connects $w$ to the root; we shall write $w\prec \ell$ if
$w\preceq v$, where  $v$ is the unique node that the line $\ell$ exits.
For any node $v$ denote by $p_{v}$ the number of lines entering $v$.

Given a rooted tree $\theta$ we denote by $N(\theta)$ the set of nodes,
by $E(\theta)$ the set of \textit{end nodes}, i.e. nodes $v$ with $p_{v}=0$,
by $V(\theta)$ the set of \textit{internal nodes}, i.e. nodes $v$ with $p_{v}\ge 1$,
and by $L(\theta)$ the set of lines; by definition $N(\theta)=E(\theta) \amalg V(\theta)$.
If, for any discrete set $A$, we denote by $|A|$ its cardinality,
we define the \textit{order} of $\theta$ as $k(\theta):=|N(\theta)|$.

We associate with each end node $v\in E(\theta)$ a \textit{mode} label $\nn_{v}\in\ZZZ^{d}$.
We split $E(\theta)=E_{0}(\theta) \amalg E_{1}(\theta)$, with
$E_{0}(\theta)=\{ v \in E(\theta) : \nn_{v} = \vzero \}$ and
$E_{1}(\theta)=\{ v \in E(\theta) : \nn_{v} \neq \vzero \}$.
With each line $\ell\in L(\theta)$ we associate
a \textit{momentum} $\nn_{\ell} \in \ZZZ^{d}$ with the constraint
\begin{equation} \nonumber \label{eq:2.9}
\nn_{\ell}=\sum_{\substack{w\in E(\theta) \\ w \prec \ell}} \nn_{w} .
\end{equation}
Finally we impose the constraints that
\begin{itemize}
\itemsep0em
\item $p_{v}\ge 2$ $\forall v\in V(\theta)$,
\item $\nn_{\ell} \neq \vzero$ for any line $\ell$ exiting a node in $V(\theta)$.
\end{itemize}

We call \textit{equivalent} two labelled rooted trees which can be transformed into
each other by continuously deforming the lines in such a way that
they do not cross each other. In the following we shall consider only
inequivalent labelled rooted trees, and we shall call them trees \textit{tout court},
for simplicity.

We associate with each node $v\in N(\theta)$  a \textit{node factor}
\begin{equation} \nonumber 
F_{v} := \begin{cases}
- \e \, a_{p_{v}} , & v \in V(\theta) , \\
\e \, f_{\nn_{v}} , & v \in E_{1}(\theta) , \\
\zeta , & v \in E_{0}(\theta) ,
\end{cases}
\end{equation}
and with each line $\ell\in L(\theta)$ a \textit{propagator}
\begin{equation} \nonumber 
\calG_{\ell} := \begin{cases}
1/D(\e, \oo\cdot\nn_{\ell}) , & \nn_{\ell} \neq \vzero , \\
1, & \nn_{\ell}=\vzero .
\end{cases}
\end{equation}
Finally we define the value of the tree $\theta$ as
\begin{equation} \label{eq:2.12}
\Val(\theta) := \Biggl( \prod_{v\in N(\theta)} F_{v} \Biggr) \Biggl( \prod_{\ell\in L(\theta)} \calG_{\ell} \Biggr) .
\end{equation}
It is not difficult to show that, with the notations above, the equations \eqref{eq:2.8}
are solved for all $k\in\NNN$, provided the coefficients $u^{(k)}_{\nn}$ are defined as
\begin{equation} \label{eq:2.13}
u^{(k)}_{\nn} = \sum_{\theta\in \TT_{k,\nn}} \Val(\theta) , \qquad \nn \in\ZZZ^{d}_{*} ,
\end{equation}
where $\TT_{k,\nn}$ is the set of non-equivalent trees of order $k$
and momentum $\nn$ associated with the root line.

For $\e$ small enough and all $s\in\RRR$, one has (see Lemma 2.2 in \cite{CFG2} for a proof)
\begin{equation} \label{eq:2.14}
|D(\e,s)| \ge \max\{|a\e|,|s|\} .
\end{equation}
Moreover, by Lemma 2.3 in \cite{CFG2}, for any tree $\theta$
one has $|E(\theta)|\ge |V(\theta)|+1$ and, as a consequence,
\begin{equation} \label{eq:2.15}
|E(\theta)|\ge \frac{1}{2} \left( k(\theta)+1 \right) .
\end{equation}
Finally set
\begin{equation} \label{eq:2.11}
C_0 := \rho^{-1}\max\{\Gamma/|a|,\Phi,1\} ,
\end{equation}
with $\rho$, $\Phi$ and $\Gamma$ defined as in \eqref{eq:2.1}.

\begin{lemma} \label{lem:2.1}
For any fixed $A\in(0,C_0)$ there exist $\bar{\e}>0$ and $\bar{\zeta}>0$
such that for any $k\ge 1$, any $\nn\in\ZZZ^{d}$ and any tree $\theta\in\TT_{k,\nn}$ one has
\begin{equation} \nonumber
\left| \Val(\theta) \right| \le A_0 \, A^{k} \prod_{v\in E(\theta)} {\rm e}^{-3\xi |\nn_{v}|/4}
\end{equation}
with $A_0$ a suitable positive constant, provided $|\e|<\bar{\e}$ and $|\zeta|<\bar{\zeta}$.
\end{lemma}

\prova
One bounds (\ref{eq:2.12}) as
\begin{eqnarray} \nonumber
\left| \Val(\theta) \right| 
& \!\!\! \le \!\!\! &
\Biggl( \prod_{v\in V(\theta)} | \e \, a_{p_{v}} | \Biggr) 
\Biggl( \prod_{v \in E_1(\theta)}  \frac{| \e \, f_{\nn_{v}} | }{|D(\oo\cdot\nn_{v},\e)|} \Biggr)
\Biggl( \prod_{v\in E_{0}(\theta)} |\zeta| \Biggr)
\Biggl( \prod_{v\in V(\theta)} \frac{1}{|a\e|} \Biggr) \nonumber \\
& \!\!\! \le \!\!\! &
|\zeta|^{|E_{0}(\theta)|} \Gamma^{|V(\theta)|} \rho^{-(|N(\theta)|-1)}
\Phi^{|E_{1}(\theta)|} |a|^{-|V(\theta)|}  
\Biggl( \prod_{v \in E_{1}(\theta)} 
\frac{ |\e| \, {\rm e}^{-\xi |\nn_{v}|} }{|D(\oo\cdot\nn_{v},\e)|} \Biggr) , \nonumber
\end{eqnarray}
where we have bounded $f_{\nn_{v}}$ as in (\ref{eq:2.1a}) and used the bound $|D(\e,s)|\ge |a\e|$
for the propagators of the lines exiting the nodes $v \in V(\theta)$.
For each end node $v \in E_{1}(\theta)$ we extract a factor ${\rm e}^{-3\xi|\nn_{v}|/4}$,
so that, if we define $C_0$ as in \eqref{eq:2.11}, we obtain
\begin{equation} \label{eq:2.16}
\left| \Val(\theta) \right| \le \rho \, C_{0}^k |\zeta|^{|E_{0}(\theta)|}
\Biggl( \prod_{v\in E_{1}(\theta)}{\rm e}^{-3\xi |\nn_{v}|/4} \Biggr)
\Biggl( \prod_{v \in E_{1}(\theta)}  \frac{ |\e| {\rm e}^{-\xi |\nn_{v}|/4} }{|D(\oo\cdot\nn_{v},\e)|} \Biggr) .
\end{equation}

For any given $n_0\in\NNN$ we have $|\oo\cdot\nn| \ge \al_{n_0}(\oo)$ for all $\nn\in\ZZZ^d_*$
such that $|\nn|\le 2^{n_0}$. Set $\de=\de(n_0):={\rm e}^{-\xi 2^{n_0}/4}$.
Let $A$ be such that $0<A<C_0$. We first fix $n_0$ such that $C_0^2 \de/|a| \le A^2$, then we fix
$\bar{\e}$ and $\bar{\zeta}$ by requiring that $C_0^2 \bar{\e}/\al_{n_0}(\oo)\le A^2$
and $C_0^2\bar{\zeta}<A^2$.

By \eqref{eq:2.14}, in \eqref{eq:2.16}, for all $v\in E_{1}(\theta)$,
we can bound  $|D(\oo\cdot\nn_{v},\e)| \ge \al_{n_0}(\oo)$ if $|\nn_{v}|\le 2^{n_0}$ and
$|D(\oo\cdot\nn_{v},\e)| \ge |\e a|$ if $|\nn_{v}|> 2^{n_0}$.  Thus,
for all $v \in E_{1}(\theta)$, one has
\begin{equation} \label{eq:2.17}
\frac{ |\e| {\rm e}^{-\xi |\nn_{v}|/4} }{|D(\oo\cdot\nn_{v},\e)|} \le
\max\left\{ \frac{\de}{|a|} , \frac{|\e|}{\al_{n_0}(\oo)} \right\} ,
\end{equation}
so that in \eqref{eq:2.16}, if $\de/|a| \ge |\e|/\al_{n_0}(\oo)$, we can bound
\begin{equation} \nonumber
C_{0}^k |\zeta|^{|E_{0}(\theta)|} \Biggl( \prod_{v \in E_{1}(\theta)}  
\frac{ |\e| \, {\rm e}^{-\xi |\nn_{v}|/4} }{|D(\oo\cdot\nn_{v},\e)|} \Biggr) \le
C_{0}^{k} \bar{\zeta}^{|E_{0}(\theta)|} \left( \frac{\de}{|a|} \right)^{|E_{1}(\theta)| } \le \frac{A}{C_0} A^{k} ,
\end{equation}
while, if $\de/|\e a| <1/\al_{n_0}(\oo)$, we can bound
\begin{equation} \nonumber 
C_{0}^k |\zeta|^{|E_{0}(\theta)|}
\Biggl( \prod_{v \in E(\theta)}  \frac{ |\e| \, {\rm e}^{-\xi |\nn_{v}|/4} }{|D(\oo\cdot\nn_{v},\e)|} \Biggr)
\le C_{0}^{k} \bar{\zeta}^{|E_{0}(\theta)|}
\left( \frac{|\e|}{\al_{n_0}(\oo)} \right)^{|E_{1}(\theta)| } \le \frac{A}{C_0} A^{k} ,
\end{equation}
where we have used twice \eqref{eq:2.15} with $k(\theta)=k$.

Summarising, for all $k\in\NNN$ and all $\nn\in\ZZZ^{d}_{*}$ we have obtained
\begin{equation} \nonumber
\left| \Val(\theta) \right| \le A_{0} 
A^{k} \Biggl( \prod_{v\in E(\theta)}{\rm e}^{-3\xi |\nn_{v}|/4} \Biggr) ,
\qquad A_{0} := \rho \frac{A}{C_{0}} .
\end{equation}
Therefore the assertion follows.
\EP

\begin{lemma} \label{lem:2.2}
For any $k\ge 1$ and $\nn\in\ZZZ^{d}_{*}$ one has
\begin{equation} \nonumber
\left| u^{(k)}_{\nn} \right| \le A_{0} C^{k} {\rm e}^{-\xi |\nn|/2} ,
\end{equation}
where $\xi$ is as in (\ref{eq:2.1a}) and $C$ is a positive constant proportional to $A$,
with $A_0$ and $A$ as in Lemma \ref{lem:2.1}.
\end{lemma}

\prova
To bound the coefficients $u^{(k)}_{\nn}$ defined by (\ref{eq:2.12})
we use the bounds of Lemma \ref{lem:2.1} and sum over all trees in $\TT_{k,\nn}$.
The sum over  the Fourier labels $\{\nn_{v}\}_{v\in E_{1}(\theta)}$ is performed by using
the factors ${\rm e}^{-3\xi |\nn_{v}|/4}$ associated with the end nodes in $E_{1}(\theta)$,
and gives a bound $C_{1}^{|E_{1}(\theta)|}{\rm e}^{-\xi |\nn|/2}$, for some positive constant $C_{1}$. 
The sum over the other labels produces a factor $C_{2}^{|N(\theta)|}$,
with $C_{2}$ a suitable positive constant.
By taking $C=AC_{1}C_{2}$ the assertion follows.
\EP

\begin{lemma} \label{lem:2.3}
For any $\oo\in\RRR^{d}$ there exist $\bar{\e}>0$ and $\bar{\zeta}>0$ such that,
for $\mu=1$, $|\e|<\bar{\e}$ and $|\zeta|<\bar{\zeta}$ the series \eqref{eq:2.7}
converges to a function $u(\pps,\e,\zeta)=u(\pps,\e,\zeta,1)$, which is analytic in $\pps$ in
a strip $\Sigma_{\x'}$, with $\xi'<\xi/2$, and such that $u(\oo t,\e,\zeta)$ solves \eqref{eq:2.5a}.
\end{lemma}

\prova
In Lemma \ref{lem:2.1} we can fix $A$ in such a way that $C\le B$, for some constant
$B<1$. Then the series \eqref{eq:2.7} converges provided $B\mu<1$, which allow $\mu=1$.
Furthermore, the function (\ref{eq:2.7}) solves (\ref{eq:2.6}) order by order by construction.
Since the series converges uniformly, then it is also a solution \textit{tout court}
to \eqref{eq:2.6} with $\mu=1$ and hence of (\ref{eq:2.5a}).  Analyticity in $\pps\in\Sigma_{\xi'}$
for any $\xi'<\xi/2$ follows from the bound on the Fourier coefficients
given by Lemma \ref{lem:2.2}.
\EP

\subsection{The bifurcation equation}

Continuity of the function $\e \mapsto u(\pps,\e,\zeta)$ holds trivially for $\e>0$.
On the contrary, continuity at $\e=0$ requires some discussion. Indeed, that $u(\pps,\e,\zeta)$ tends to $0$
as $\e\to 0$ does not follow from Lemma \ref{lem:2.3}, since the constants $A$ and $A_0$
do not tend to $0$ as $\e\to0$. However, continuity at $\e=0$ can be proved
by following the same lines as in the proof of Lemma \ref{lem:2.2}, up to some minor changes.

\begin{lemma} \label{lem:2.4}
Let $\bar{\e}$ and $\bar{\zeta}$ be as in Lemma \ref{lem:2.3}. For any $|\zeta|<\bar{\zeta}$,
the function $u(\pps,\e,\zeta)$ in Lemma \ref{eq:2.3} is continuous in $\e\in[0,\bar{\e})$.
In particular, it tends to $0$ as $\e\to0$.
\end{lemma}

\prova
Let $\zeta$ be such that $|\zeta|<\bar{\zeta}$.
As already noted, continuity in $\e$ is obvious for $\e>0$.
Set
\begin{equation} \nonumber
F(\e,\zeta) = \|u(\cdot,\e,\zeta)\|_{\io} := \sup\{ u(\pps,\e,\zeta) : \pps \in \Sigma_{\xi'} \} ,
\end{equation}
with $\xi'$ as in Lemma \ref{lem:2.3}. Since $F(0,\zeta)=0$,
we have only to prove that $F(\e,\zeta)  \to 0$ as $\e\to 0$,
that is that for all $\eta>0$ there exists $\de>0$ such that $0<\e<\de$ implies $|F(\e,\zeta)|<\eta$.

Let $\bar{\e}$ be as in Lemma \ref{lem:2.3}. The series in \eqref{eq:2.7} with $\mu=1$ can be written as
\begin{equation} \nonumber
u(\pps, \e, \zeta) = \sum_{\nn \in \ZZZ^{d}} {\rm e}^{i\nn\cdot\pps} u_{\nn} ,
\qquad u_{\nn} := \sum_{k=1}^{\io} u^{(k)}_{\nn} ,
\end{equation}
so that
\begin{equation} \nonumber
F(\e,\zeta) \le \sum_{k=1}^{\io} \sum_{\nn\in\ZZZ^{d}} \bigl| u^{(k)}_{\nn} \bigr| \, {\rm e}^{\x'|\nn|} .
\end{equation}
By reasoning as in the proof of Lemma \ref{lem:2.1} -- see in particular \eqref{eq:2.16} --, we can bound
\begin{eqnarray} 
\sum_{\nn\in\ZZZ^{d}} \bigl| u^{(k)}_{\nn} \bigr| \, {\rm e}^{\x'|\nn|}
& \!\!\! \le \!\!\! &
\sum_{\nn\in\ZZZ^{d}} \sum_{\theta \in \TT_{\nn,k}} \bigl| \Val(\theta) \bigr| \, {\rm e}^{\x'|\nn|} \le
\rho \, C_{0}^{k} \sum_{\theta \in \TT_{\nn,k}} \!\! |\zeta|^{|E_{0}(\theta)|} \!\!\!\! \prod_{v \in E_{1}(\theta)}
\sum_{\nn_v \in \ZZZ^{d}} \frac{|\e| \, {\rm e}^{-\xi |\nn_v|/4}}{|D(\oo\cdot\nn_v,\e)|} \nonumber \\
& \!\!\! \le \!\!\! &
\rho \, C_{0} A^{-1} C^{k} \!\! \sum_{\nn \in \ZZZ^{d}} \frac{|\e|\,{\rm e}^{-\xi |\nn|/4}}{|D(\oo\cdot\nn,\e)|} , \nonumber
\end{eqnarray}
where we have used the fact that $|E_{1}(\theta)|\ge 1$ (since $\nn\neq\vzero$) and
the bound \eqref{eq:2.17} for all the end nodes $v\in E_{1}(\theta)$ but one.
Therefore we obtain, for any $N\in\NNN$,
\begin{equation} \nonumber \label{eq:2.18}
F(\e,\zeta) \le \frac{\rho \, C_0 \,A^{-1}}{1-C} \sum_{\substack{\nn\in\ZZZ^d \\ |\nn| \le N}}
\frac{|\e|\,{\rm e}^{-\xi |\nn|/4}}{|D(\oo\cdot\nn,\e)|} + D_0 {\rm e}^{-\xi N/8} , \qquad D_0 :=
\frac{\rho \, C_0 \, A^{-1}}{1-C} \sum_{\nn \in \ZZZ^{d}} {\rm e}^{-\xi|\nn|/8} .
\end{equation}

Fix $\eta>0$. Choose $N$ such that $D_0 {\rm e}^{-\xi N/8} < \eta/2$. If we define
\begin{equation} \nonumber
r_{N} : = \min\{ |\oo\cdot\nn| : 0<|\nn| \le N \} ,
\end{equation}
%
we can bound
\begin{equation} \nonumber
\frac{\rho \, C_0 \, A^{-1}}{1-C} \sum_{\substack{\nn\in\ZZZ^d \\ |\nn| \le N}}
\frac{|\e|\,{\rm e}^{-\xi |\nn|/4}}{|D(\oo\cdot\nn,\e)|} \le \frac{|\e|\, D_{0} }{r_{N}} .
\end{equation}
Thus, for $\de$ small enough and $0<\e<\de$,
we have $|\e|D_{0}/r_{N} < \eta /2$  and hence $|F(\e)| < \eta$.
\EP

Note that we are not able to prove more than continuity for the function $\e \mapsto u(\pps,\e,\zeta)$.
Indeed, $\partial_{\e} u(\pps,\e,\zeta)$ is well defined for $\e>0$, but the argument used in proving
Lemma \ref{lem:2.4} does not allow us to obtain its boundedness  as $\e\to 0$.

\begin{lemma} \label{lem:2.5}
Let $\bar{\e}$ and $\bar{\zeta}$ be as in Lemma \ref{lem:2.2} and let $u=u(\oo t,\e,\zeta)$
be as in Lemma \ref{lem:2.3}. There exist neighbourhoods $U \subset(-\bar{\e},\bar{\e})$
and $V \subset (-\bar{\zeta},\bar{\zeta})$ and a function $\zeta : U \to V$ such that for all $\e\in U$
the equation \eqref{eq:2.5b} holds for $\zeta=\zeta(\e)$. Moreover the function $\e \mapsto \zeta(\e)$
is continuous in $U$ and $\zeta(\e)$ is the only solution to \eqref{eq:2.5b} in $V$.
\end{lemma}

\prova
Write \eqref{eq:2.5b} as
\begin{equation} \nonumber
H(\zeta,\e) := a \, \zeta + \left[ G(c_0+\zeta+u) \right]_{\vzero} = 0 .
\end{equation}
By construction, the function $u(\pps,\e,\zeta)$ is analytic on $\zeta$
in a neighbourhood of the origin. Therefore,
$H(\zeta,\e)$ is analytic in $\zeta$ and, by Lemma \ref{lem:2.3},  is continuous in $\e$. One has
\begin{equation} \nonumber
H(0,0) = 0 , \qquad \frac{\partial}{\partial \zeta} H(0,0) = a \neq 0 ,
\end{equation}
so that we can apply the implicit function theorem, in the version of
Loomis and Sternberg \cite{LS},  so as to conclude that there exist neighbourhoods 
$U\subset(-\bar{\e},\bar{\e})$ and $V\subset(-\bar{\zeta},\bar{\zeta})$
such that for all $\e\in U$ one can find a unique value $\zeta(\e) \in V$,
depending continuously on $\e$, such that $H(\zeta(\e),\e)=0$.
\EP

\begin{lemma} \label{lem:2.6}
Let $u(\pps,\e,\zeta)$ and $\zeta(\e)$ be as in Lemma \ref{lem:2.3}
and in Lemma \ref{lem:2.5}, respectively.
There exists $\e_{0}>0$ such that for all $\e\in(0,\e_0)$ the function 
$x(t,\e)=c_0+\zeta(\e)+u(\oo t,\e,\zeta(\e))$ solves \eqref{eq:2.2}.
Moreover $x(t,\e) \to c_0$ as $\e\to 0$.
\end{lemma}

\prova
The result follows immediately from Lemma \ref{lem:2.3} and Lemma \ref{lem:2.5}.
\EP

\zerarcounters 
\section{Proof of Theorem \ref{thm:2}}
\label{sec:3} 

Let the function $h(\pps,x)$ in \eqref{eq:1.2} be analytic on the domain
$\Sigma_{\xi} \times \Delta(c_0,\rho_0)$, with the same notations as in Section \ref{sec:1}.
We expand
\begin{equation} \nonumber
h(\pps,x) = \sum_{\nn\in\ZZZ^{d}} h_{\nn}(x) \, {\rm e}^{i\nn\cdot\pps} , \qquad
h_{\nn}(x) = \sum_{p=0}^{\io} a_{\nn,p} \,  \left(x - c_{0} \right)^{p} , \qquad 
a_{\nn,p} = \frac{1}{p!} \frac{\partial^p h_{\nn}}{\partial x^{p}} (c_0) , 
\end{equation}
so that, for any $\rho<\rho_0$,
\begin{equation} \label{eq:3.0}
\left| a_{\nn,p} \right| \le \Gamma \rho^{-p} {\rm e}^{-\xi|\nn|} ,
\end{equation}
for a suitable positive constant $\Gamma$.
Then we rewrite \eqref{eq:1.2} as
\begin{eqnarray}
\e \ddot x + \dot x + \e \, a \left( x - c_0 \right) 
& \!\!\! + \!\!\! & 
\e \sum_{\nn\in\ZZZ^{d}_{*}} {\rm e}^{i\nn\cdot\oo t} h_{\nn}(c_0) \nonumber \\
& \!\!\! + \!\!\! &  
\e \sum_{\nn\in\ZZZ^{d}_{*}} a_{\nn,1} \, {\rm e}^{i\nn\cdot\oo t} \left( x - c_0 \right) +
\e \sum_{p=2}^{\io} \sum_{\nn\in\ZZZ^{d}} a_{\nn,p} \, {\rm e}^{i\nn\cdot\oo t} 
\left( x - c_0 \right)^{p} = 0,\nonumber 
\end{eqnarray}
where we have used that $h_{\vzero}(c_0)=0$ and $a:=a_{\vzero,1} \neq 0$
by Hypothesis \ref{hyp:2} with $\gotn=1$.

We look for a solution of the form \eqref{eq:2.3}. By passing to Fourier space, setting
\begin{equation} \nonumber
\al_{1}(\pps) := \sum_{\nn\in\ZZZ^{d}_{*}} a_{\nn,1} \, {\rm e}^{i\nn\cdot\pps} , \qquad
\al_{p}(\pps) := \sum_{\nn\in\ZZZ^{d}} a_{\nn,p} \, {\rm e}^{i\nn\cdot\pps} ,
\end{equation}
and defining  $D(\e,s)$ as in \eqref{eq:2.4}, we obtain the equations
\begin{subequations} \label{eq:3.1}
\begin{align}
D(\e,\oo\cdot\nn) \, u_{\nn} & = - \e h_{\nn}(c_0) - 
\e \left[ \al_{1} \left( x-c_0 \right) \right]_{\nn} - \e \sum_{p=2}^{\io}
\left[ \al_{p} \left( x-c_0 \right)^{p} \right]_{\nn} , \qquad \nn \neq \vzero ,
\label{eq:3.1a} \\
\e \, a \, \zeta  & = - \left[ \al_{1} \left( x-c_0 \right) \right]_{\vzero} - \e \sum_{p=2}^{\io}
\left[ \al_{p} \left( x-c_0 \right)^{p} \right]_{\vzero} .
\label{eq:3.1b}
\end{align}
\end{subequations}

By following the same strategy as in Section \ref{sec:2}, we shall study first
\eqref{eq:3.1a} and look for a solution depending on the parameter $\zeta$.
Thereafter we shall fix $\zeta$ by requiring that \eqref{eq:3.1b} is solved as well.


Instead of \eqref{eq:3.1a} we consider the equation
\begin{equation} \label{eq:3.2}
D(\e,\oo\cdot\nn) \, u_{\nn} = - \mu \, \e h_{\nn}(c_0) -  \mu \, \e
\left[ \al_{1} \left( x-c_0 \right) \right]_{\nn} - \mu \, \e\sum_{p=2}^{\io}
\left[ \al_{p} \left( x-c_0 \right)^{p} \right]_{\nn} , \qquad \nn \neq \vzero ,
\end{equation}
and look for a quasi-periodic solution to (\ref{eq:3.2}) in the form of a power series in $\mu$,
\begin{equation} \label{eq:3.3}
x(t,\e,\mu) = c_{0} + \zeta + u(\oo t, \e,\zeta,\mu) , \qquad
u(\oo t, \e,\zeta,\mu) = \sum_{k=1}^{\io} \sum_{\nn\in\ZZZ^{d}_{*}} 
\mu^{k} {\rm e}^{\ii\nn\cdot\pps} u^{(k)}_{\nn} .
\end{equation}
We find the recursive equations
\begin{subequations} \label{eq:3.4}
\begin{align}
D(\e, \oo\cdot\nn) \, u^{(1)}_{\nn} & = - \e \, h_{\nn}(c_0) 
\label{eq:3.4a} \\
D(\e, \oo\cdot\nn) \, u^{(k)}_{\nn} & = - \e \, \sum_{\nn_0 \in\ZZZ^{d}_{*}}
a_{\nn_0,1} \, u^{(k-1)}_{\nn-\nn_0} 
\nonumber \\
&  - \e \, \sum_{p=2}^{\io} \!\!\!\!\!\!
\sum_{\substack{ k_{1},\ldots,k_{p} \ge 1 \\ k_{1}+\ldots+k_{p}=k-1}}
\sum_{\substack{\nn_{0},\nn_{1},\ldots,\nn_{p} \in\ZZZ^{d} \\ \nn_{0}+\nn_{1}+\ldots+\nn_{p}=\nn}} \!\!\!\!\!\!
a_{\nn_{0},p} \, u^{(k_{1})}_{\nn_{1}} \ldots u^{(k_{p})}_{\nn_{p}} , \qquad k \ge 2 ,
\label{eq:3.4b}
\end{align}
\end{subequations}
where $\nn\neq\vzero$ and we have set once more $u^{(1)}_{\vzero}=\zeta$ and $u^{(k)}_{\vzero}=0$ $\forall k\ge 2$.

We have still a tree representation of the coefficients $u^{(k)}_{\nn}$, with a few differences with respect to Section \ref{sec:2}.
Define the sets $N(\theta)$, $E(\theta)$, $E_{0}(\theta)$, $E_{1}(\theta)$, $V(\theta)$ and $L(\theta)$ as previously
and call $k(\theta):=|N(\theta)|$ the order of $\theta$.
Now we split $V(\theta)=V_{1}(\theta) \amalg V_{2}(\theta)$, with
$V_{1}(\theta)=\{ v \in V(\theta) : p_{v} = 1 \}$ and $V_{2}(\theta)=\{ v \in E(\theta) : p_{v} \geq 2 \}$.
Contrary to Section \ref{sec:2}, now in general $V_{1}(\theta) \neq \emptyset$.

We associate with each node $v\in N(\theta)$ a \textit{mode} label $\nn_{v}\in\ZZZ^{d}$
and with each line $\ell\in L(\theta)$ a \textit{momentum} $\nn_{\ell} \in \ZZZ^{d}$ with the constraint
\begin{equation} \nonumber  
\nn_{\ell}=\sum_{\substack{w\in N(\theta) \\ w \prec \ell}} \nn_{w} .
\end{equation}
Finally we impose the constraints that
\begin{itemize}
\itemsep0em
\item $\nn_{v} \neq \vzero$ $\forall v\in V_{1}(\theta)$,
\item $\nn_{\ell} \neq \vzero$ for any line $\ell$ exiting a node in $V(\theta)$.
\end{itemize}

We associate with each node $v\in N(\theta)$  a \textit{node factor}
\begin{equation} \label{3.6}
F_{v} := \begin{cases}
- \e \, a_{\nn_{v},p_{v}} , & v \in V(\theta) , \\
-\e \, h_{\nn_{v}}(c_0) , & v \in E_{1}(\theta) , \\
\zeta , & v \in E_{0}(\theta) ,
\end{cases}
\end{equation}
and with each line $\ell\in L(\theta)$ a \textit{propagator}
\begin{equation} \label{eq:3.7}
\calG_{\ell} := \begin{cases}
1/D(\e, \oo\cdot\nn_{\ell}) , & \nn_{\ell} \neq \vzero , \\
1, & \nn_{\ell}=\vzero .
\end{cases}
\end{equation}
We define the value of the tree $\theta$ as \eqref{eq:2.12} and write
$u^{(k)}_{\nn}$ as in \eqref{eq:2.13}, where $\TT_{k,\nn}$ denotes the set of non-equivalent trees
of order $k$ and momentum $\nn$ associated with the root line, constructed according to the
new rules. Then the coefficients $u^{(k)}_{\nn}$ solve formally \eqref{eq:3.4}.

In a tree $\theta$ we define a \emph{chain} $\CCCC$ as a subset of $\theta$ formed
by a maximal connected set of nodes $v \in V(\theta)$ with $p_{v}=1$ and by the lines exiting them.
Therefore, if $V(\CCCC)$ and $L(\CCCC)$ denote the set of nodes and the set of lines of $\CCCC$ and 
$V(\CCCC)=\{v_1,v_2,\ldots,v_p\}$, with $v_1 \succ v_2 \succ \ldots \succ v_p$,
then $L(\CCCC)=\{\ell_{1},\ell_{2},\ldots,\ell_{p}\}$, where
$\ell_{i}$ is the line exiting $v_{i}$, for $i=1,\ldots,p$, and $\ell_{i}$ enters the node $v_{i-1}$ for $i=2,\ldots,p$.
We call $p=|V(\CCCC)|=|L(\CCCC)|$ the \emph{length} of the chain $\CCCC$. Finally we define
the \emph{value} of the chain $\CCCC$ as
\begin{equation} \nonumber
\Val(\CCCC) := \Biggl( \prod_{v\in V(\CCCC)} F_{v} \Biggr) 
\Biggl( \prod_{\ell\in L(\CCCC)} \calG_{\ell} \Biggr) 
\end{equation}
and denote by $\gotC(\theta)$ the set of all the chains contained in $\theta$.
The main difference with respect to the trees considered in Section \ref{sec:2} is that, now,
the trees may contain chains.

We define $n_{0}$ and $\de$ as in Section \ref{sec:2}. Define also
$\be:=\max\{\de,2|\e a|/\al_{n_0}(\oo)\}$. Finally set
\begin{equation} \label{eq:3.11}
C_{0}=\rho^{-1}\max\{\Gamma/|a|,1\} ,
\end{equation}
with $\rho$ and $\Gamma$ defines as in \eqref{eq:3.0}.

\begin{lemma} \label{lem:3.1}
Let $\CCCC$ be a chain of length $p\ge 1$. Then
\begin{equation} \nonumber
|\Val(\CCCC)| \le C_{0}^{p} \be^{(p-1)/2} \prod_{v\in V(\CCCC)} {\rm e}^{-3\xi|\nn_{v}|/4} .
\end{equation}
\end{lemma}

\prova
We proceed by induction on $p$. 
If $p=1$ then $\CCCC$ contains only one node $v$, so that
\begin{equation} \nonumber
\left| \Val(\CCCC) \right| \le \frac{|\e| \Gamma 
\rho^{-1} {\rm e}^{-\xi|\nn_{v}|}}{|D(\e,\oo\cdot\nn_{v})|} \le C_{0} {\rm e}^{-3\xi|\nn_{v}|/4} ,
\end{equation}
where we have bounded $|D(\e,\oo\cdot\nn_{v})| \ge |\e a|$ by \eqref{eq:2.14}.

If $p\ge 2$, let  $v_1 \succ v_2 \succ v_3 \succ \ldots \succ v_{p}$ be the nodes
in $\CCCC$ and let $\ell_1,\ell_2,\ell_3,\ldots,\ell_p$ be the lines exiting such nodes.
The nodes $\{v_2,v_3,\ldots,v_p\}$ and the lines $\{\ell_2,\ell_3,\ldots,\ell_p\}$
form a chain $\CCCC'$ of length $p-1$ and, if $p\ge 3$, the nodes $\{v_3,\ldots,v_p\}$
and the lines $\{\ell_3,\ldots,\ell_p\}$ form a chain $\CCCC''$ of length $p-2$.

We assume that the bound holds up to $p-1$.
If $|\oo\cdot\nn_{\ell_1}| \ge \al_{n_0}(\oo)/2$, then one has
\begin{equation} \nonumber
\left| \Val(\CCCC) \right| \le \frac{|\e| \Gamma 
\rho^{-1} {\rm e}^{-\xi|\nn_{v_1}|}}{|D(\e,\oo\cdot\nn_{\ell_1})|} 
\left| \Val(\CCCC') \right| \le C_{0} \be \left( C_{0}^{p-1} \be^{(p-2)/2} \right) 
\prod_{v\in V(\CCCC)} {\rm e}^{-3\xi|\nn_{v}|/4} ,
\end{equation}
which yields the bound for $p$.
If $|\oo\cdot\nn_{\ell_1}| < \al_{n_0}(\oo)/2$ and $p\ge 3$ we distinguish between two cases.
If $|\oo\cdot\nn_{\ell_2}| \ge \al_{n_0}(\oo)/2$, we can bound
\begin{equation} \nonumber
\left| \Val(\CCCC) \right| \le 
\frac{|\e| \Gamma \rho^{-1} {\rm e}^{-\xi|\nn_{v_1}|}}{|D(\e,\oo\cdot\nn_{\ell_1})|} 
\frac{|\e| \Gamma \rho^{-1} {\rm e}^{-\xi|\nn_{v_2}|}}{|D(\e,\oo\cdot\nn_{\ell_2})|} 
\left| \Val(\CCCC'') \right| \le C^2_{0} \be \left( C_{0}^{p-2} \be^{(p-3)/2} \right) 
\prod_{v\in V(\CCCC)} {\rm e}^{-3\xi|\nn_{v}|/4} ,
\end{equation}
so that the bound follows once more.

If $|\oo\cdot\nn_{\ell_2}| < \al_{n_0}(\oo)/2$, then
\begin{eqnarray}
|\oo\cdot\nn_{v_1}| & \!\!\! = \!\!\! & |\oo\cdot\nn_{v_1}+\oo\cdot\nn_{\ell_2} - \oo\cdot\nn_{\ell_2}|
\le |\oo\cdot\nn_{v_1}+\oo\cdot\nn_{\ell_2} | + |\oo\cdot\nn_{\ell_2}| \nonumber \\
& \!\!\! = \!\!\! & |\oo\cdot\nn_{\ell_1} | + |\oo\cdot\nn_{\ell_2}| < \al_{n_0}(\oo) , \nonumber
\end{eqnarray}
so that, since $\nn_{v_1} \neq \vzero$ and hence $\oo\cdot\nn_{v_{1}} \neq 0$,
we conclude that $|\nn_{v_1}| > 2^{n_0}$, which allows us to bound
${\rm e}^{-\xi|\nn_{v_1}|} \le \de {\rm e}^{-3\xi|\nn_{v_1}|/4}$. Therefore we obtain
\begin{equation} \nonumber
\left| \Val(\CCCC) \right| \le 
C_{0}^{2} \de \left| \Val(\CCCC'') \right| \le C^2_{0} \de \left( C_{0}^{p-2} \be^{(p-3)/2} \right) 
\prod_{v\in V(\CCCC)} {\rm e}^{-3\xi|\nn_{v}|/4} ,
\end{equation}
which gives the bound for $p$ in this case too.

Finally if $p=2$ we can reason as in the case $p\ge 3$, the only difference being that
the quantity $|\Val(\CCCC'')|$ has to replaced with $1$ -- since there is no further
chain $\CCCC''$ for $p=2$. Therefore we obtain
$|\Val(\CCCC)| \le C_{0}^2 \be$, which is the desired bound.
\EP

\begin{lemma} \label{lem:3.2}
For any tree $\theta$ one has
\begin{enumerate}
\itemsep0em
\item \label{1}$|E(\theta)| \ge |V_{2}(\theta)|+1$,
\item \label{2}$|\gotC(\theta)| \le |E(\theta)| + |V_{2}(\theta)|$,
\item \label{3}$|\gotC(\theta)| \le |V_{1}(\theta)|$,
\item \label{4}$k(\theta)=|E(\theta)|+|V_{1}(\theta)|+|V_{2}(\theta)|$.
\end{enumerate}
\end{lemma}

\prova
Property \ref{1} can be proved as the inequality $|E(\theta)| \ge |V(\theta)|+1$
in Section \ref{sec:2}. Property \ref{2} is easily proved by noting that each chain
has to be preceded by either an end node or a node $v\in V_{2}(\theta)$.
Property \ref{3} is trivial, since any chain has to contain at least one node $v\in V_{1}(\theta)$.
Finally property \ref{4} follows from the definition of order. 
\EP

A result analogous to Lemma \ref{lem:2.1} still holds. This can be proved as follows.
For any tree $\theta$ we have
\begin{eqnarray} \nonumber
\left| \Val(\theta) \right| 
& \!\!\! \le \!\!\! &
\Biggl( \prod_{v\in V_{2}(\theta)} \frac{| \e \, a_{p_{v},\nn_{v}}|}{|\e a|} \Biggr) 
\Biggl( \prod_{v \in E_1(\theta)}  \frac{| \e \, h_{\nn_{v}}(c_0) | }{|D(\oo\cdot\nn_{v},\e)|} \Biggr)
\Biggl( \prod_{v\in E_{0}(\theta)} |\zeta| \Biggr)
\Biggl( \prod_{\CCCC \in \gotC(\theta)} |\Val(\CCCC)| \Biggr) \nonumber \\
& \!\!\! \le \!\!\! &
\rho \, C_{0}^{k} |\zeta|^{|E_{0}(\theta)|} 
\Biggl( \prod_{v \in E_{1}(\theta)} 
\frac{ |\e| \, {\rm e}^{-\xi |\nn_{v}|/4} }{|D(\oo\cdot\nn_{v},\e)|} \Biggr) 
\Biggl( \prod_{\CCCC \in \gotC(\theta)} \be^{(|V(\CCCC)|-1)/2} \Biggr) 
\Biggl( \prod_{v \in V(\theta)} {\rm e}^{-3\xi |\nn_{v}|/4} \Biggr), \nonumber
\end{eqnarray}
so that
\begin{equation} \nonumber
\left| \Val(\theta) \right| \le \rho \, C_{0}^{k} \bar{\zeta}^{|E_{0}(\theta)|} 
\left( \max \left\{ \frac{\de}{|a|},\frac{\bar{\e}}{\al_{n_0(\oo)}} \right\} \right)^{|E_{1}(\theta)|}
\be^{(|V_1(\theta)|-|\gotC(\theta)|)/2} 
\Biggl( \prod_{v \in V(\theta)} {\rm e}^{-3\xi |\nn_{v}|/4} \Biggr)
\end{equation}
Therefore, for any constant $A\in(0,C_0)$, if we fix first $n_0$ and hence $\bar{\e}$ and $\bar{\zeta}$
so that for any $|\e|<\bar{\e}$ and any $|\zeta|<\bar{\zeta}$ one has
\begin{equation} \label{eq:3.8}
C_{0}^{4} \max \left\{ \bar{\zeta}, \frac{\de}{|a|}, \frac{\bar{\e}}{\al_{n_0(\oo)}} ,\beta \right\} < A^4 ,
\end{equation}
we obtain
\begin{equation} \label{eq:3.9}
\left| \Val(\theta) \right| \le \rho \, C_{0}^{k}
\left( \frac{A}{C_{0} }\right)^{4|E(\theta)| + 2|V_{1}(\theta)| - 2 |\gotC(\theta)|}
\Biggl( \prod_{v \in V(\theta)} {\rm e}^{-3\xi |\nn_{v}|/4} \Biggr) .
\end{equation}
By using Lemma \ref{lem:3.2} we can bound
\begin{eqnarray} 
4|E(\theta)| + 2|V_{1}(\theta)| - 2 |\gotC(\theta)|
& \!\!\! = \!\!\! & 
2|E(\theta)| + 2|E(\theta)| + |V_{1}(\theta)| - |\gotC(\theta)| + \left( |V_{1}(\theta)| - |\gotC(\theta)| \right) \nonumber \\
& \!\!\! \ge \!\!\! & 
2|E(\theta)| + 2 \left( |V_{2}(\theta)| + 1 \right)  + |V_{1}(\theta)| - |\gotC(\theta)|
\nonumber \\
& \!\!\! \ge \!\!\! & 
|E(\theta)| + |V_{2}(\theta)| + |V_{1}(\theta)| + 2 +  \left( |E(\theta)| + |V_{2}(\theta)| - |\gotC(\theta)|
\right) \nonumber \\
& \!\!\! \ge \!\!\! & 
|E(\theta)| + |V_{2}(\theta)| + |V_{1}(\theta)| + 2 = k(\theta) + 2 , \nonumber
\end{eqnarray}
which, inserted into \eqref{eq:3.9}, gives
\begin{equation} \nonumber 
\left| \Val(\theta) \right| \le \rho \, C_{0}^{k}
\left( \frac{A}{C_{0} }\right)^{k+2} \Biggl( \prod_{v \in V(\theta)} {\rm e}^{-3\xi |\nn_{v}|/4} \Biggr) \le A_{0}  A^{k} 
\Biggl( \prod_{v \in V(\theta)} {\rm e}^{-3\xi |\nn_{v}|/4} \Biggr) ,
\qquad A_{0} := \rho \frac{A^2}{C_{0}^2} .
\end{equation}
From here on, we can reason as in Section \ref{sec:2} and we find that the coefficients
$u^{(k)}_{\nn}$ can be bounded as
\begin{equation} \nonumber
\left| u^{(k)}_{\nn} \right| \le A_{0} C^{k} {\rm e}^{-\xi |\nn|/2} ,
\end{equation}
for a suitable constant $C$ proportional to the constant $A$ in \eqref{eq:3.8}.
Thus, by choosing $\bar{\e}$ and $\bar{\zeta}$ small enough, we can make $C$
such that  $C \le B <1$, so that the series \eqref{eq:3.3} converges for $\mu=1$.
Therefore the function $x(t,\e,\zeta)$ solves the range equation \eqref{eq:3.1a}
for any $\zeta$ small enough.

Both the proof of continuity in $\e$ of the function $x(t,\e,\zeta)$
and the discussion of the bifurcation equation in order to fix the parameter $\zeta$
can be repeated exactly as in Section \ref{sec:2}, so we omit the details.

\appendix



\end{document}